\documentclass{article}

\def\a{\alpha}
\def\b{\beta}

\def\G{\Gamma}
\def\D{\triangle}
\def\t{\tau}
\def\d{\delta}
\def\th{\theta}
\def\l{\lambda}
\def\L{\Lambda}

\def\s{\sigma}

\def\f{\rightarrow}
\def\q{\forall}

\def\v{\vdash}

\def\p{\succ}

\def\ou{\vee}
\def\et{\wedge}

\newcommand{\cqfd}{\hbox{}\nobreak\hfill\nobreak$\spadesuit$}

\begin{document}

\begin{center}
{\Large\bf 
LES TYPES DE DONN\'EES SYNTAXIQUES DU SYST\`EME ${\cal F}$}\\[0.3cm]

{\Large\bf  SYNTACTICAL DATA TYPES OF SYSTEM ${\cal F}$}\\[0.3cm]

{\large par Samir FARKH et Karim NOUR}\footnote{LAMA - \'Equipe de Logique, Universit\'e de Savoie -
73376 Le Bourget du Lac cedex - Email knour@univ-savoie.fr} 
\end{center}

R\'esum\'e. - {\it Nous pr\'esentons dans ce papier une d\'efinition
purement syntaxique des types entr\'ees et des types sorties du
syst\`eme ${\cal F}$. Nous d\'efinissons les types de donn\'ees
syntaxiques comme \'etant des types entr\'ees et sorties. Nous
d\'emontrons que les types \`a quantificateurs positifs sont des types
de donn\'ees syntaxiques et qu'un type entr\'ee est un type
sortie. Nous imposons des restrictions sur la r\`egle d'\'elimination
des quantificateurs pour d\'emontrer qu'un type sortie est un type entr\'ee.}  \\[0.3cm]

Abstract. - {\it We give in this paper a purely syntactical definition
of input and output types of system ${\cal F}$. We define the
syntactical data types as input and output types. We show that any
type with positive quantifiers is a syntactical data type and that an input
type is an output type. We give some restrictions on the
$\q$-elimination rule in order to prove that an output type is an
input type.}  \\



{\Large\bf Introduction} \\

Le syst\`eme de typage ${\cal F}$ a \'et\'e introduit par J.-Y. Girard
(voir [3]).  Ce syst\`eme est bas\'e sur le calcul propositionnel
intuitionniste du second ordre, et donc donne la possibilit\'e de
quantifier sur les types. En plus du th\'eor\`eme de normalisation
forte qui assure la terminaison des programmes, le syst\`eme ${\cal
F}$ permet d'une part, d'\'ecrire des programmes pour toutes les
fonctions dont la terminaison est d\'emontrable dans l'arithm\'etique
de Peano du second ordre, et d'autre part, de d\'efinir tous les types
de donn\'ees courants : bool\'eens, entiers, listes d'objet, etc.\\

Nous avons essay\'e de trouver une d\'efinition syntaxique \`a un type
de donn\'ees dans le syst\`eme ${\cal F}$. Nous l'avons d\'efini comme
\'etant un type {\bf entr\'ee} et {\bf sortie}.
  
\begin{itemize}
\item{\bf Type entr\'ee :} Il faut qu'une machine soit capable de
tester si les entr\'ees sont bien typ\'ees, c'est \`a dire que le
probl\`eme de typage d'un type entr\'ee soit d\'ecidable.  On
d\'efinit donc un type entr\'ee comme \'etant un type dont toutes les
d\'emonstrations se font dans une restriction d\'ecidable (not\'ee
${\cal F}_0$) du syst\`eme ${\cal F}$.
 
\item{\bf Type sortie :} Si $E$, $S$ sont deux types du syst\`eme
${\cal F}$, et $t$ un $\l$-terme clos tel que $\v_{\cal F} t : E \f
S$, alors $t$ peut \^etre vu comme un programme qui \`a un \'el\'ement
de $\L (E)$ (l'ensemble des $\l$-termes de type E) associe un
\'el\'ement de $\L (S)$. Pour nous, la sortie doit d\'ependre de
l'entr\'ee.  Les seuls programmes (fonctions) qui ne tiennent pas
compte de leurs arguments (entr\'ees) sont les fonctions constantes
(i.e les $\l$-termes de la forme $\l xt$ o\`u $t$ est un terme clos).
En formalisant les types qui v\'erifient cette propri\'et\'e, nous
avons obtenu la d\'efinition suivante : un type {\it sortie} est un
type clos $S$ tel que si $\v \l xt : \q X (X\f S)$ ($t$ est un
$\l$-terme normal), alors $x$ est non libre dans $t$.
\end{itemize}

Nous avons remarqu\'e que les types de donn\'ees ainsi d\'efinis
contiennent les types \`a quantificateurs positifs. Ensuite, nous
avons montr\'e qu'un type entr\'ee est un type sortie. Ce r\'esultat
nous a conduit \`a regarder la r\'eciproque.  Nous l'avons
d\'emontr\'e dans des cas particuliers, o\`u on impose des
restrictions sur la r\`egle d'\'elimination des quantificateurs du
second ordre. Le cas g\'en\'eral, qui reste une conjecture, est
d\'emontr\'e lorsqu'on se restreint aux termes du $\l I$-calcul (voir
[9]). Enfin nous avons obtenu un r\'esultat sur les op\'erateurs de
mise en m\'emoire de J.-L. Krivine (voir [6]).  Nous avons montr\'e
que si $D$ est un type de donn\'ees syntaxique, alors un $\l$-terme de
type $D^* \f \neg \neg D$ ($D^*$ est la traduction de G\"odel de $D$)
est un op\'erateur de mise en m\'emoire pour $D$.  Le fait que $D$ est
un type sortie est n\'ecessaire pour avoir ce r\'esultat.
    
\section {Notations et d\'efinitions}

On d\'esignera par $\L$ l'ensemble des termes du $\l$-calcul pur, dits
aussi {\bf $\l$-termes}.  \'Etant donn\'es des $\l$-termes $t, u,
u_1,..., u_n$, l'application de $t$ \`a $u$ sera not\'ee $(t)u$, et
$(... ((t)u_1) ...)u_n$ sera not\'e $(t)u_1...u_n$.  Si $t$ est un
$\l$-terme, on d\'esigne par $Fv (t)$ l'ensemble de ses variables
libres. On note par $\f_{\b}$ la $\b$-r\'eduction. Un $\l$-terme $t$
soit poss\`ede un {\bf redex de t\^ete faible} [i.e. $t = (\l x u) v
v_1 ... v_m$, le redex de t\^ete faible est $(\l x u ) v$], soit est
en {\bf forme normale de t\^ete faible} [i.e. $t= (x ) v_1 ... v_m$ ou
$t = \l x v$]. La notation $u \p_f v$ signifie que $v$ est obtenu \`a
partir de $u$ par {\bf r\'eduction de t\^ete faible}.\\

Les {\bf types} du syst\`eme ${\cal F}$ sont les formules construites
\`a l'aide d'un ensemble d\'enombrable de variables propositionnelles
$X, Y$,..., et deux connecteurs $\f$ et $\q$. \'Etant donn\'es un
$\l$-terme $t$, un type $A$, et un contexte $\G = \{x_1 : A_1,..., x_n
: A_n \}$, on d\'efinit au moyen des r\`egles suivantes la notion
``$t$ est typable de type $A$ dans le contexte $\G$''. Cette notion
est not\'ee $\G \v_{\cal F} t : A$.
\begin{itemize}
\item $(ax)$ $\G \v_{\cal F} x_i : A_i$  $ (1 \leq i \leq n)$.
\item $(\f_i)$ Si $\G, x : B \v_{\cal F} t : C$, alors $\G \v_{\cal F} \l xt : B \f C$.
\item$(\f_e)$ Si $\G \v_{\cal F} u : B\f C$, et $\G \v_{\cal F} v : B$, alors
$\G \v_{\cal F} (u)v : C$.
\item $(\q_i)$ Si $\G \v_{\cal F} t : A$, et $X$ ne figure pas dans $\G$,
alors $\G \v_{\cal F} t : \q XA$.
\item $(\q_e)$ Si $\G \v_{\cal F} t : \q XA$, alors, pour tout type $C$, $\G \v_{\cal F} t : A[C/X]$.
\end{itemize}

Le syst\`eme ${\cal F}$ poss\`ede les propri\'et\'es suivantes (voir
[4]) : \\

{\bf Proposition 1.1} {\it (i) Si $\G \v_{\cal F} t : A$, et
$t\f_\beta t'$, alors $\G \v_{\cal F} t' : A$.

(ii) Si $\G \v_{\cal F} t : A$, alors $t$ est fortement
normalisable.}\\

{\bf Proposition 1.2} {\it (i) Si $\G \v_{\cal F} t : A$, alors pour
tout type $G$, $\G [G/X] \v_{\cal F} t : A [G/X]$.

(ii) Si $\G, x : B \v_{\cal F} u : A$ et $\G \v_{\cal F} v : B$, alors
$\G \v_{\cal F} u [v/x] : A$}.\\

On ne consid\`ere dans ce papier que des types ``propres'' (c'est \`a
dire les variables sur lesquelles on quantifie figurent dans le
type).\\

Dans la suite, on note par $\q {\bf X} A$ la formule $\q X_1 ... \q
X_n A$ ($n \geq 0$).\\
  
Une partie $G$ de $\L$ est dite {\bf satur\'ee} si, quels que soient
les $\l$-termes $t$ et $u$, on a : ($u \in G$ et $t \p_f u) \Rightarrow t
\in G$. Il est clair que l'intersection d'un ensemble de parties
satur\'ees de $\L$ est satur\'ee. \'Etant donn\'ees deux parties $G$
et $G'$ de $\L$, on d\'efinit une partie de $\L$, not\'ee $G\f G'$, en
posant : $ u\in (G\f G') \Leftrightarrow (u)t \in G'$ quel que soit
$t\in G$. Si $G'$ est satur\'ee, alors $G\f G'$ est satur\'ee pour
toute partie $G \subset \L$. Une {\bf interpr\'etation} $I$ est, par
d\'efinition, une application $X \f |X|_I$ de l'ensemble des variables
de type dans l'ensemble des parties satur\'ees de $\L$. $X$ \'etant
une variable de type, et $G$ une partie satur\'ee de $\L$, on
d\'efinit une interpr\'etation $J = I [X \leftarrow G]$ en posant
$|X|_J = G$, et $|Y|_J = |Y|_I$ pour toute variable $Y \neq X$. Pour
chaque type $A$, sa valeur $|A|_I$ dans l'interpr\'etation $I$ est une
partie satur\'ee d\'efinie comme suit, par induction sur $A$ :

- si $A$ est une variable de type, $|A|_I$ est d\'ej\`a d\'efinie ;

- $|A \f B|_I = |A|_I \f |B|_I$ ;

- $|\q XA|_I = \cap \{ |A|_{I [X \leftarrow G]}$ pour toute partie
satur\'ee $G \}$.\\

Pour tout type $A$, on note $|A| = \cap \{|A|_I$ ; $I$
interpr\'etation$\}$. \\

Le r\'esultat suivant est connu sous le nom du lemme d'ad\'equation
(voir [5]).\\

{\bf Th\'eor\`eme 1.3} {\it Soient $A$ un type et $t$
un $\l$-terme clos. Si $\v_{\cal F} t : A$, alors $t \in |A|$.}\\

Les types \`a {\bf quantificateurs positifs} (resp. \`a {\bf
quantificateurs n\'egatifs)}, not\'es en abr\'eg\'e $\q^+$ (resp.
$\q^-$), sont d\'efinis de la mani\`ere suivante : 

- une variable propositionnelle $X$ est $\q^+$ et $\q^-$ 

- si $A$ est $\q^+$ (resp. $\q^-$) et $B$ est $\q^-$ (resp. $\q^+$),
alors $B \rightarrow A$ est $\q^+$ (resp. $\q^-$) 

- si $A$ est $\q^+$ et $X$ est libre dans $A$, alors $\q XA$ est $\q^+$\\

Le r\'esultat suivant constitue une sorte de r\'eciproque du
th\'eor\`eme 1.3 (voir [2]). \\

{\bf Th\'eor\`eme 1.4} {\it Soient A un type $\q^+$ du syst\`eme
${\cal F}$, et $t$ un $\l$-terme, alors $t \in |A|$ ssi il existe un
$\l$-terme clos $t'$ tel que $t \f_{\b} t'$ et $\v_{\cal F} t' :
A$.}\\

Le syst\`eme de typage {\bf simple} ${\cal S}$ est la restriction du
syst\`eme ${\cal F}$ aux types qui ne contiennent pas de
quantificateurs. Ce syst\`eme poss\`ede donc trois r\`egles de typage
: $(ax)$, $(\f_i)$ et $(\f_e)$. Si $t$ est un $\l$-terme, $A$ un type,
et $\G = \{x_1 : A_1,..., x_n : A_n \}$ un contexte, alors on \'ecrit
$\G \v_{\cal S} t : A$ ssi $t$ est typable dans le syst\`eme ${\cal
S}$ de type $A$ \`a partir de $\G$.

\section{Types de donn\'ees syntaxiques} 

\subsection {Types sorties}
     
{\bf D\'efinition :} Un type clos $S$ est un {\bf type sortie} ssi
pour tout $\l$-terme normal $t$, si $\v_{\cal F} \l xt : \q X (X \f
S)$, alors  $x \not \in Fv (t)$.\\
 
Cela veut dire que les fonctions (programmes) \`a valeurs dans un type
sortie $S$ ind\'ependamment du type de leurs arguments (entr\'ees)
sont les fonctions constantes.\\
   
Soit $O$ une constante de type. La d\'efinition d'un type sortie est
\'equivalente \`a la suivante :\\
 
{\bf D\'efinition :} Un type clos $S$ (ne contenant pas la constante
$O$) est un type sortie ssi pour tout $\l$-terme normal $t$, 
si $\a : O \v_{\cal F} t : S$, alors  $\a \not \in Fv (t)$. \\
 
{\bf Exemples :}
\begin{itemize}
\item[(1)] Les types $Id = \q X (X \f X)$ (type de l'identit\'e), $B
= \q X \{X\f (X\f X)\}$ (type des bool\'eens), et $N = \q X \{X\f
[(X\f X)\f X] \}$ (type des entiers) sont des types sorties.

On va faire la preuve pour le type $B$. Soit $t$ un $\l$-terme normal
tel que $\a : O \v_{\cal F} t : B$, donc $\a : O \v_{\cal F} t : O'\f
(O' \f O')$, o\`u $O'$ est une constante de type diff\'erente de
$O$. $t$ est donc une abstraction, soit $t = \l xu$, d'o\`u $\a : O, x
: O' \v_{\cal F} u :O'\f O'$. $u$ s'\'ecrit n\'ecessairement $\l yv$,
avec $\a : O, x : O', y : O' \v_{\cal F} v : O'$, donc $v = x$ ou $v =
y$, ce qui fait que $t = \l x\l yx =$ {\bf 1} ou $t = \l x\l yy =$
{\bf 0}.  

\item[(2)]  Le type $D = \q X \{\q Y (Y\f X) \f X \}$ n'est pas un
type sortie. En effet, posons $t = \l x (x)\a$; $t$ est un $\l$-terme
normal non clos. D'autre part on a $\a : O, x : \q Y (Y\f X) \v_{\cal
F} x : \q Y (Y\f X)$, donc $\a : O, x : \q Y (Y\f X) \v_{\cal F} x :
O\f X$.  D'o\`u $\a : O, x : \q Y (Y\f X) \v_{\cal F} (x)\a : X$, et
donc $\a : O \v_{\cal F} t = \l x (x)\a : D$.

\item[(3)] Le type $N\f N$ (type des fonctions d'entiers dans les
entiers) n'est pas un type sortie. En effet, il suffit de trouver un
$\l$-terme $t_{\a}$ normal non clos tel que $\a : O \v_{\cal F} t_{\a}
: N\f N$.  On a, $\G = x : X, \a : O, n : N, z : X\f X \v_{\cal F} n :
N$, donc $\G \v_{\cal F} n : (O\f X) \f [((O\f X) \f (O\f X))\f (O\f
X)]$, par suite $\G \v_{\cal F} (n)\l yx : [((O\f X) \f (O\f X))\f
(O\f X)]$. Par cons\'equent $\G \v_{\cal F} (((n)\l yx)\l xx)\a : X$, et
donc $\a : O \v_{\cal F} \l n\l x\l z (((n)\l yx)\l xx)\a : N \f N$.

\item[(4)] Soit $S$ un type du syst\`eme ${\cal F}$. Si l'ensemble
des $\l$-termes normaux de type $S$ est fini, alors $S$ est un type
sortie. En effet, si $S$ n'est pas un type sortie, alors il existe un
$\l$-terme normal $t$ contenant $\a$ tel que $\a : O \v_{\cal F} t :
S$. D'o\`u, d'apr\`es la proposition 1.2, $\a : O [E/O] \v_{\cal F} t
: S$ pour tout type $E$. Donc si $u$ est un $\l$-terme clos tel que
$\v_{\cal F} u : E$, alors, d'apr\`es la proposition 1.2, $\v_{\cal F}
t [u/ \a] : S$. Or $\a$ ne peut pas \^etre en position d'application,
car $\a$ est de type atomique.  Donc, comme $t$ est normal, alors $t
[u/ \a]$ est normal.  On obtient donc un nombre infini de $\l$-termes
normaux de type $S$.

\item[(5)] Si $E \f F$ est un type sortie, alors $F$ est un type
sortie. En effet, soit $t$ un $\l$-terme normal tel que $\a : O
\v_{\cal F} t : F$, alors $\a : O, y : E \v_{\cal F} t : F$, en
choisissant une variable $y$ non libre dans $t$. Donc $\a : O \v_{\cal
F} \l yt : E \f F$. Comme $E\f F$ est un type sortie, alors $\a \not
\in Fv (\l yt)$, et donc $\a \not \in Fv (t)$.  
\end{itemize}
  
{\bf D\'efinition :} Soit $K$ un variable ou une constante de type. On
  dit qu'un type {\bf $A$ se termine par $K$} ssi $A$ est obtenu par
  les r\`egles suivantes :

 - $K$ se termine par $K$ ;

 - si $A$ se termine par $K$, alors $B \f A$ se termine par $K$, pour
  tout type $B$ ;

 - si $A$ se termine par $K$, alors $\q XA$ se termine par $K$, pour
  toute variable de type $X \not = K$.\\

Un type $A$ qui se termine par $K$ s'\'ecrit alors : $A = \q
\mbox{\boldmath$X_0$} (B_1 \f \q \mbox{\boldmath$X_1$} (B_2 \f \q
\mbox{\boldmath$X_2$} (...(B_r \f \q \mbox{\boldmath$X_r$} (B_{r+1}\f
K))...)))$.\\
 
 On se propose de d\'emontrer le th\'eor\`eme suivant :\\

{\bf Th\'eor\`eme 2.1.1} {\it Un type clos $S$ est un type sortie ssi
 pour tout $\l$-terme normal $t$ et pour tous types $A_1$,..., $A_r$
 qui se terminent par $O$, si $x_1 : A_1,..., x_r : A_r \v_{\cal F} t :
 S$, alors $ x_i \not \in Fv (t)$ pour tout $1 \leq i \leq r$}. \\

 On a besoin du lemme suivant :\\

{\bf Lemme 2.1.2} {\it Soient $t$ un $\l$-terme normal, $v$ un
 $\l$-terme, et $\a, x$ deux variables tels que $x \in Fv (t)$. Si $t
 [\l x_1...\l x_n\a /x] \f_{\b} v$, alors $\a \in Fv (v)$.}\\

{\bf Preuve :} Par induction sur $t$.

-- Si $t = x$, alors $t [\l x_1...\l x_n\a /x] = \l x_1...\l x_n\a =
v$. Donc $\a \in Fv (v)$.

-- Si $t = (x)u_1...u_m$, $u_i$ normal, $1 \leq i \leq m$, alors $t
[\l x_1...\l x_n\a /x] =$

$ (\l x_1...\l x_n\a)u_1...u_m$.

- Si $m = n$, alors $t [\l x_1...\l x_n\a /x] \f_{\b} \a$.

- Si $m > n$, alors $t [\l x_1...\l x_n\a /x] \f_{\b}
(\a)u_{n+1}...u_m$.

- Si $m < n$, alors $t [\l x_1...\l x_n\a /x] \f_{\b} \l x_{m+1}...\l
x_n \a$.

 Dans les trois cas on a bien $\a \in Fv (v)$.

 -- Si $t = (y)u_1...u_m$, $u_i$ normal, et $y$ une variable
diff\'erent de $x$.  Comme $x \in Fv (t)$, alors il existe $i$, $1
\leq i \leq m$ tel que $x \in Fv (u_i)$. Donc $t [\l x_1...\l x_n\a
/x] = (y)u'_1...u_i [\l x_1...\l x_n\a /x]...u'_m \f_{\b} \\
(y)u'_1...v'...u'_m$, avec $u_i [\l x_1...\l x_n\a /x] \f_{\b}
v'$. D'o\`u, d'apr\`es l'hypoth\`ese d'induction, $\a \in Fv (v')$, et
par cons\'equent $\a \in Fv (v)$.

 -- Si $t = \l zu$, alors $ \l zu [\l x_1...\l x_n\a /x] \f_{\b} v$,
donc $v = \l zw$, o\`u $u [\l x_1...\l x_n\a /x]$\\ $\f_{\b} w$. Or $x \in
Fv (t)$, donc $x \in Fv (u)$, et d'apr\`es l'hypoth\`ese d'induction,
$\a \in Fv (w)$, d'o\`u $\a \in Fv (v)$. \cqfd\\

On a donc imm\'ediatement :\\
 
{\bf Corollaire 2.1.3} {\it Soient $t$ un $\l$-terme normal, $v$ un
$\l$-terme clos, et $\a, x$ deux variables. Si $t [\l x_1...\l x_n\a
/x] \f_{\b} v$, alors $x \not \in Fv (t)$}.\\

Plus g\'en\'eralement, on a le r\'esultat suivant : \\

{\bf Lemme 2.1.4} {\it Soient $t$ un $\l$-terme normal, et $v$ un
$\l$-terme clos.\\ Si $t [\l y_{11}...\l y_{1n_1}\a /x_1,..., \l
y_{r1}...\l y_{rn_r}\a /x_r ] \f_{\b} v$, alors $x_i \not \in Fv (t)$
pour tout $1\leq i \leq r$.}\\

On peut alors d\'eduire la preuve du th\'eor\`eme 2.1.1. \\

{\bf Preuve du th\'eor\`eme 2.1.1 :} La condition suffisante est
\'evidente. Supposons que $x_1 : A_1,..., x_r : A_r \v_{\cal F} t :
S$, avec $S$ un type sortie et $t$ un $\l$-terme normal.  Comme $A_i$
se termine par $O$, alors $A_i$ s'\'ecrit :\\ $A_i = \q
\mbox{\boldmath$X_1$} (B_1 \f \q \mbox{\boldmath$X_2$} (B_2 \f \q
\mbox{\boldmath$X_3$} (...(B_{n_i-2} \f \q \mbox{\boldmath$X_{n_i-1}$}
(B_{n_i-1}\f O))...)))$, donc $\a : O \v_{\cal F} \l y_{i1}...\l
y_{in_i} \a : A_i$. \\D'o\`u, $\a : O \v_{\cal F} t [\l y_{11}...\l
y_{1n_1}\a /x_1,..., \l y_{r1}...\l y_{rn_r}\a /x_r] : S$, d'apr\`es
la proposition 1.2. Par cons\'equent, $t [\l y_{11}...\l y_{1n_1}\a
/x_1,..., \l y_{r1}...\l y_{rn_r}\a /x_r] \f_{\b} v$, avec $v$ un
$\l$-terme clos. D'o\`u, d'apr\`es le lemme 2.1.4, $x_i \not \in Fv
(t)$ pour tout $1\leq i \leq r$. \cqfd\\

{\bf D\'efinition :} Si $A$ et $B$ sont deux types du syst\`eme ${\cal
F}$, alors le type $A \et B = \q X \{ (A \f (B\f X))\f X \}$ est dit
{\bf type produit} de $A$ et $B$, le type $A \ou B = \q X \{ (A\f X)
\f ( (B\f X) \f X)\}$ est dit {\bf type somme disjointe} de $A$ et
$B$, et le type $LA = \q X\{X\f [(A \f (X\f X))\f X]\}$ est dit {\bf
type liste d'objet de} $A$.\\
 
{\bf Corollaire 2.1.5} {\it Si $A$ et $B$ sont des types sorties,
alors $A \et B$, $A \ou B$, et $LA$ sont des types sorties.}\\
 
{\bf Preuve :} Faisons la preuve pour $A \et B$. Soit $t$ un
 $\l$-terme normal tel que $\a : O \v_{\cal F} t : A\et B$, donc $\a :
 O \v_{\cal F} t : (A, B\f O)\f O$, d'o\`u $t = \l xu$, avec $\a : O,
 x : (A, B\f O) \v_{\cal F} u : O$. Comme $O$ est une constante de
 type, on obtient $u = (x)ab$, avec $\a : O, x : (A, B\f O) \v_{\cal
 F} a : A$ et $\a : O, x : (A, B\f O) \v_{\cal F} b : B$. Or $A$ et
 $B$ sont des types sorties, d'o\`u d'apr\`es le th\'eor\`eme 2.1.1,
 $a$ et $b$ sont des $\l$-termes clos, et donc $t$ est clos. \cqfd\\
    
On va d\'emontrer que tout type $\q^+$ du syst\`eme ${\cal F}$ est un
type sortie. On a besoin du lemme suivant d\'emontr\'e dans [7].\\

{\bf Lemme 2.1.6} {\it Soient $t$ un $\l$-terme normal, $A_1, ...,
 A_n$ des types $\q^-$, $S$ un type $\q^+$, $O$ une constante de type
 qui ne figure pas dans $A_1, ..., A_n, S$, et $B_1, ..., B_m$ des
 types qui se terminent par $O$.  Si $x_1 : A_1, ..., x_n : A_n, y_1 :
 B_1, ..., y_m : B_m \v_{\cal F} t : S$, alors $y_i \not \in Fv (t)$,
 pour tout $1\leq i \leq m$}.\\

{\bf Th\'eor\`eme 2.1.7} {\it Si $S$ est un type $\q^+$, alors $S$ est
un type sortie.}\\

{\bf Preuve :} Il suffit d'appliquer le lemme 2.1.6, avec $n = 0, m =
1$ et $B_1 = O$. \cqfd \\

{\bf Remarque :} La r\'eciproque du th\'eor\`eme 2.1.7 n'est pas en
g\'en\'eral vraie. En effet, consid\'erons le type $S = \q X \{\q Y
(Y\f X) \f Id\}$. $S$ n'est pas $\q^+$, d'autre part si $t$ est un
$\l$-terme normal tel que $\a : O \v_{\cal F} t : S$, alors $t = \l
xu$, avec $\a : O, x : \q Y (Y\f O) \v_{\cal F} u : Id$.  Comme $Id$
est un type sortie, $\a \not \in Fv (u)$, d'o\`u $\a \not \in Fv (t)$,
et par cons\'equent $S$ est un type sortie.\\
   
{\bf D\'efinition\footnote{Cette d\'efinition est due \`a
J.-L. Krivine (voir [4])}:} Un type clos $A$ du syst\`eme ${\cal F}$
est un {\bf type de donn\'ees} ssi $|A| \neq \emptyset$ et tout terme
$t \in |A|$ se r\'eduit par $\b$-r\'eduction \`a un terme clos.\\

{\bf Exemples :} J.-L. Krivine a montr\'e dans [4] que les types $Id$, $B$
et $N$ sont des types de donn\'ees.\\

{\bf Th\'eor\`eme 2.1.8} {\it Si $A$ est un type de donn\'ees,
alors $A$ est un type sortie.}\\

{\bf Preuve :} Soit $t$ un $\l$-terme normal tel que $\a : O \v_{\cal
F} t : A$. On d\'efinit une interpr\'etation ${\cal I}$ en posant
$|O|_{\cal I} = \{\t \in \L : \t \p_f \a \}$.  On a $\a \in |O|_{\cal
I}$, donc, d'apr\`es le lemme d'ad\'equation, $t \in |A|_{\cal I} =
|A|$. Comme $A$ est un type de donn\'ees, alors $t$ est un terme clos,
ce qui fait que $A$ est sortie. \cqfd \\

{\bf Th\'eor\`eme 2.1.9} {\it Si $A$ est un type clos, $\q^+$ et
d\'emontrable, alors $A$ est un type de donn\'ees.}\\
 
{\bf Preuve :} Comme $A$ est d\'emontrable, alors il existe un
$\l$-terme clos $t$, tel que $ \v_{\cal F} t : A$. D'o\`u d'apr\`es le
lemme d'ad\'equation $t \in |A|$, et donc $|A| \neq \emptyset$. Soit
donc $t \in |A|$. Comme $A$ est un type $\q^+$, alors d'apr\`es le
th\'eor\`eme 1.4, il existe un $\l$-terme clos $t'$ tel que $t \f_{\b}
t'$ et $ \v_{\cal F} t' : A$.  Par cons\'equent $A$ est un type de
donn\'ees. \cqfd \\
   
{\bf Remarque :} Il existe des types de donn\'ees qui ne sont pas
    $\q^+$.  Consid\'erons par exemple le type $S = \q X \{\q Y (Y\f
    X)\f Id \}$.  Il est clair que $S$ n'est pas $\q^+$. De plus $S$
    est un type de donn\'ees. En effet, si $\v_{\cal F} t : S$, alors
    $\v_{\cal F} t : \q Y (Y\f X) \f Id$, et donc $t = \l xu$, avec $x
    : \q Y (Y\f X) \v_{\cal F} u : Id$. D'o\`u $u = \l yv$, avec $x :
    \q Y (Y\f X), y : Z \v_{\cal F} v : Z$, $Z$ \'etant une variable
    de type.  Par cons\'equent $v = y$ et $t = \l x\l yy$. 

 D'autre part, soit $t \in |S|$ ; $x$ \'etant une variable du
    $\l$-calcul qui n'est pas libre dans $t$, on d\'efinit une
    interpr\'etation ${\cal I}$ en posant $|X|_{\cal I} = \{\t \in \L$
    ; il existe $G$ une partie satur\'ee et $u \in G$ tels que $\t
    \f_{\b} (x)u \}$.  $|X|_{\cal I}$ est \'evidemment une partie
    satur\'ee. On a $t \in |\q Y (Y\f X) \f Id|_{\cal I}$. Or $x \in
    |(Y\f X)|_{J = I [Y \leftarrow G]}$ pour toute partie satur\'ee
    $G$ de $\L$, car si $v \in G$, alors $(x)v \in | X |_{\cal I}$,
    par d\'efinition de ${\cal I}$. Donc $(t)x \in |Id| = \{\t \in \L;
    \t \f_{\b} \l yy \}$, d'o\`u $(t)x \f_{\b} \l yy$, et par
    cons\'equent $(t)x$ est normalisable et donc $t$ est
    normalisable. Soit $t'$ sa forme normale, alors deux cas peuvent
    se pr\'esenter : 

 -- Si $t'$ commence par $\l$, on \'ecrit $t' = \l xu$, donc $(t)x
    \f_{\b} (t')x \f_{\b} u$, d'o\`u $u = \l yy$. On a donc $t\f_{\b}
    t' = \l x \l yy$.

 -- Sinon $t' = (f)u_1...u_n$, donc $(t)x \f_{\b} (t')x =
    (f)u_1...u_nx$.  D'o\`u $(f)u_1...u_nx = \l yy$, ce qui est
    impossible. \\ On vient donc de d\'emontrer que $t \in |S|$ ssi $t
    \f_{\b} \l x\l yy$ et $\v_{\cal F} \l x\l yy : S$.
 
\subsection{Types entr\'ees}

{\bf D\'efinition :} On d\'efinit le syst\`eme ${\cal F}_0$ comme
\'etant le syst\`eme ${\cal F}$ sans la r\`egle de typage $(\q_e)$.\\

{\bf Th\'eor\`eme 2.2.1} {\it Le probl\`eme de typabilit\'e d'un
$\l$-terme normal dans le syst\`eme ${\cal F}_0$ est d\'ecidable.}\\

{\bf Preuve} : Ceci provient des \'equivalences suivantes :

(i) $\G \v_{{\cal F}_0} x : A$ ssi $x : A \in \G$.

 (ii) $\G , x : B \v_{{\cal F}_0} (x)t_1 ... t_n : A$ ssi $B =
B_1,...,B_n \f A$ et $\G , x : B \v_{{\cal F}_0} t_i : B_i$ ($1 \leq i
\leq n$).

 (iii) $\G \v_{{\cal F}_0} \l x t : A$ ssi $A = \q {\bf X} (B \f C)$
et $\G , x : B \v_{{\cal F}_0} t : C$.\cqfd \\

Le r\'esultat suivant est d\'emontr\'e dans [7] : \\

{\bf Th\'eor\`eme 2.2.2} {\it Soient A un type $\q^+$ du syst\`eme
${\cal F}$, et $t$ un $\l$-terme normal clos.  Si $\v_{\cal F} t :
A$, alors $\v_{{\cal F}_0} t : A$.}\\

{\bf D\'efinition :} Un type clos $E$ du syst\`eme ${\cal F}$ est dit
{\bf entr\'ee} s'il v\'erifie la propri\'et\'e suivante : si $t$ est
un $\l$-terme normal tel que $\v_{\cal F} t : E$, alors $\v_{{\cal
F}_0} t : E$. \\
  
Cela veut dire qu'un type entr\'ee est un type dont toutes les
d\'emonstrations se font dans le syst\`eme ${{\cal F}_0}$. D'apr\`es
le th\'eor\`eme 2.2.2, on a le r\'esultat suivant : \\
 
{\bf Th\'eor\`eme 2.2.3} {\it Si $A$ est un type $\q^+$, alors $A$ est
un type entr\'ee}.\\

{\bf Exemples :} D'apr\`es le th\'eor\`eme 2.2.3, les types $Id = \q X
(X \f X)$ et $B = \q X \{X\f (X\f X)\}$ sont des types entr\'ees. Par
contre le type $D = \q X \{\q Y (Y\f X) \f X \}$ ne l'est pas, puisque
$\v_{\cal F} \l x (x)\l yy : D$ et $\not \v_{{\cal F}_0} \l x (x)\l yy : D$.\\
     
{\bf D\'efinition :} On dit qu'un type clos $D$ du syst\`eme ${\cal
F}$ est un {\bf type de donn\'ees syntaxique}, s'il est \`a la fois un
type entr\'ee et sortie.\\

D'apr\`es les th\'eor\`emes 2.1.7 et 2.2.3, on a :\\

{\bf Th\'eor\`eme 2.2.4} {\it Si $D$ est un type clos, $\q^+$ du
syst\`eme ${\cal F}$, alors $D$ est un type de donn\'ees
syntaxique.}\\

 On va montrer qu'un type entr\'ee est un type sortie, et donc les
 types de donn\'ees syntaxiques seront les types sorties. Dans la
 preuve on a besoin de trois lemmes.\\

On note par $p_n$, $n \in$ ${\bf N}$, le $\l$-terme $\l x_1...\l x_n
 \l xx$. Le lemme suivant est facile \`a d\'emontrer.\\
 
{\bf Lemme 2.2.5} {\it Soit $G$ un type du syst\`eme ${\cal S}$.  Si
 $\G \v_{\cal S} p_n : G$, alors $G = C_1\f (...\f (C_n \f D)...)$ et
 $\G, x_1 : C_1, ..., x_n : C_n \v_{\cal S} \l xx : D$.}\\
 
{\bf D\'efinition :} On d\'efinit la {\bf longueur d'un type} $E$
(qu'on note $Lg (E)$) comme \'etant le nombre des $\f$ dans $E$. \\

{\bf Lemme 2.2.6} {\it Soient $E$ un type du syst\`eme ${\cal S}$,
$A_1$,..., $A_m$, $G$ des sous-types de $E$, $t$ un $\l$-terme normal,
et $\a$ une variable libre de $t$ qui n'est pas en position
d'application.  Si $\G = x_1 : A_1,..., x_m : A_m \v_{\cal S} t [p_n /
\a] : G$, alors $Lg (E) \geq n$.}\\
  
{\bf Preuve :} Par induction sur $t$.

-- Si $t = \a$, alors $t [p_n / \a] = p_n$, par cons\'equent $\G
\v_{\cal S} p_n : G$. D'o\`u, d'apr\`es le lemme 2.2.5, $G = C_1,...,
C_n \f D$ et $\G, x_1 : C_1, ..., x_n : C_n \v_{\cal S} \l xx : D$,
donc $Lg (G) \geq n$. Or $G$ est un sous-type de $E$, donc $Lg (E)
\geq n$.

-- Si $t = \l xu$, alors $\G \v_{\cal S} \l xu [p_n / \a] : G$. Donc
$G = C\f D$ et $\G, x : C \v_{\cal S} u [p_n / \a] : D$, d'o\`u,
d'apr\`es l'hypoth\`ese d'induction, $Lg (E) \geq n$.

-- Si $t = (x_i) v_1...v_k$ ($k \geq 1$), comme $\a \in Fv (t)$, alors
il existe $j$, $1 \leq j \leq k$ tel que $\a \in Fv (v_j)$. Donc $t
[p_n / \a] = (x_i) v'_1...v_j [p_n / \a]...v'_k$, et $\G \v_{\cal S}
(x_i) v'_1...v_j [p_n / \a]...v'_k : G$, d'o\`u $A_i = C'_1\f (...\f
(C_j\f ...\f (C'_k \f G)...)...)$, avec $\G \v_{\cal S} v_j [p_n / \a]
: C_j$.  Donc $C_j$ est un sous-type de $E$, et d'apr\`es
l'hypoth\`ese d'induction, $Lg (E) \geq n$. \cqfd \\
 
{\bf D\'efinition :} Soit $A$ un type du syst\`eme ${\cal F}$.  On
d\'efinit le type $A^s$ par induction sur $A$ :

 - si $A = X$, alors $A^s = X$ ;

 - si $A = B\f C$, alors $A^s  = B^s \f C^s$ ;

 - si $A = \q XB$, alors $A^s = B^s$.\\
 
Soit $\G = x_1 : A_1, ..., x_n : A_n$, on note $\G^s = x_1 : A^s_1,
 ..., x_n : A^s_n$. Alors on a le lemme suivant : \\
  
{\bf Lemme 2.2.7} {\it Soient $E$ un type du syst\`eme ${\cal F}$, et
$t$ un $\l$-terme.  Si $\G = x_1 : A_1, ..., x_n : A_n \v_{{\cal F}_0}
t : E$, alors $\G^s \v_{\cal S} t : E^s$.}\\

{\bf Preuve :} Par induction sur le typage.  \cqfd \\
  
{\bf Th\'eor\`eme 2.2.8} {\it Si $E$ est un type entr\'ee du syst\`eme
${\cal F}$, alors $E$ est un type sortie}.\\

{\bf Preuve :} Soit $t$ un $\l$-terme normal tel que $\a : O \v_{\cal
F} t : E$.  Supposons que $Lg (E) = r$, et consid\'erons le $\l$-terme
$p_n$, avec $n > r$. Si $\a \in Fv (t)$, alors $\a : X_1,..., X_n\f Id
\v_{\cal F} t : E$.  Comme $\v_{\cal F} p_n : X_1,..., X_n\f Id$,
alors, d'apr\`es la proposition 1.2, $\v_{\cal F} t [p_n / \a] : E$.
Or $\a$ ne peut pas \^etre en position d'application, car $\a$ est de
type atomique, donc $t [p_n / \a]$ est normal. $E$ \'etant un type
entr\'ee, on obtient alors $\v_{{\cal F}0} t [p_n / \a] : E$. Donc,
d'apr\`es le lemme 2.2.7, $\v_{\cal S} t [p_n / \a] : E^s$, d'o\`u,
d'apr\`es le lemme 2.2.6, $Lg (E) \geq n$. Ce qui fait que $r \geq n$,
contradiction. Par cons\'equent $\a \not \in Fv (t)$, et donc $E$ est
un type sortie. \cqfd\\
 
{\bf Th\'eor\`eme 2.2.9} {\it Si $A$ et $B$ sont des types entr\'ees,
alors $A \et B$, $A \ou B$, et $LA$ sont des types entr\'ees.}\\
 
{\bf Preuve :} Faisons la preuve pour $A \et B$. Soit $t$ un
 $\l$-terme normal tel que $\v_{\cal F} t : A\et B$, donc $t = \l
 x(x)ab$, avec $\a : O, x : (A, B\f O) \v_{\cal F} a : A$ et $\a : O,
 x : (A, B\f O) \v_{\cal F} b : B$. Or $A$ et $B$ sont des types
 entr\'ees donc sorties. D'apr\`es le th\'eor\`eme 2.1.1, on a
 $\v_{\cal F} a : A$ et $\v_{\cal F} b : B$, donc $\v_{{\cal F}_0} a :
 A$ et $\v_{{\cal F}_0} b : B$. Par cons\'equent $\v_{{\cal F}_0} t :
 A\et B$.\cqfd
     
\section{Sortie $\Rightarrow$ Entr\'ee} 

 Nous avons montr\'e dans le paragraphe pr\'ec\'edent qu'un type
entr\'ee est un type sortie.  Ce r\'esultat nous a conduit \`a
regarder la r\'eciproque. Nous allons la d\'emontrer dans des cas
particuliers, o\`u on impose des restrictions sur la r\`egle de typage
($\q_e$). \\
  
{\bf D\'efinition :} Si $G$ est un type du syst\`eme ${\cal F}$, alors
on note par $G^o$ le type $O \f G \et O$.\\
  
{\bf Pr\'esentation globale de la preuve :}\\
Soit $E$ un type sortie.  On se propose de trouver des conditions pour
que $E$ soit un type entr\'ee.  Raisonnons par l'absurde, donc
supposons que $E$ n'est pas un type entr\'ee, alors il existe un
$\l$-terme normal $T$ tel que $\v_{\cal F} T : E$ et $\not \v_{{\cal
F}_0} T : E$.  Ce qui fait qu'au cours de typage de $T$ on a utilis\'e
au moins une fois la r\`egle de typage ($\q_e$).  C'est \`a dire, on
avait dans le typage de $T$, $\displaystyle\frac{ \G \v_{\cal F}
\d : \q XA } { \G \v_{\cal F} \d : A[G/X] }$, $\d$ \'etant un
sous-terme de $T$, $\G$ un contexte et $A$, $G$ deux types, avec $X$
figure dans $A$. Pour aboutir \`a une contradiction, il suffit de
trouver un terme normal $T^o$ contenant $\a$ tel que $\a : O \v_{\cal
F} T^o : E$. On proc\`ede de la mani\`ere suivante: On reprend le
typage de $T$ et lorsqu'on arrive \`a $\G \v_{\cal F} \d : \q XA$, on
applique la r\`egle ($\q_e$), en rempla\c cant cette fois $X$ par
$G^o$. C'est \`a dire, on obtient $\G \v_{\cal F} \d : A[G^o/X]$. Donc
s'il existe un terme $T'_A$ tel que $\a : O \v_{\cal F} T'_A : A
[G^o/X] \f A [G/X]$, on aura $\G, \a : O \v_{\cal F} (T'_A)\d :
A[G/X]$.  Puis on suivra les m\^emes \'etapes que pr\'ec\'edemment
(dans le typage de $T$) pour obtenir un terme. $T^o$ est la forme
normale du terme obtenu.\\
   
Le premier lemme \`a d\'emontrer est donc l'existence d'un terme
$T'_A$ tel que $\a : O \v_{\cal F} T'_A : A [G^o/X] \f A [G/X]$. \\

{\bf D\'efinition :} Si $F$ est un type, on d\'efinit deux $\l$-termes
$T_F$ et $T'_F$ par induction sur $F$ de la fa\c con suivante :

- si $X \not \in Fv (F)$, alors $T_F = T'_F = \l xx$ ;

- si $F = X$, alors $T_F = \l x \l {\b} \l g (g)x{\a}$ et $T'_F = \l x
(x) \a {\bf 1}$ ;

- si $F = C \f D$, alors $T_F = \l x\l y (T_D) (x) (T'_C)y$ et $T'_F =
 \l x\l y (T'_D) (x) (T_C)y$ ;

- si $F = \q YB$, alors $T_F = \l x (T_B)x$ et $T'_F = \l x (T'_B)x$.\\
   
On a alors le lemme suivant : \\
    
{\bf Lemme 3.1} {\it $\a : O \v_{\cal F} T_A : A [G/X] \f A [G^o/X]$
 et $\a : O \v_{\cal F} T'_A : A [G^o/X] \f A [G/X]$.} \\

 {\bf Preuve :} Par induction sur $A$. 

-- Si $A = X$, alors $A [G/X] = G$ et $A [G^o/X] = G^o$.  On a $\a :
O, x : G \v_{\cal F} \l g (g)x{\a} : G \et O$, donc $\a : O \v_{\cal
F} \l x \l {\b} \l g (g)x{\a} : G \f G^o$.  D'autre part, $x : G^o, \a
: O \v_{\cal F} (x)\a : G \et O$, donc $x : G^o, \a : O \v_{\cal F}
(x)\a {\bf 1} : G$. D'o\`u $\a : O \v_{\cal F} \l x(x)\a {\bf 1} :
G^o\f G$.

-- Si $A = B \f C$, on a par hypoth\`ese d'induction $\a : O \v_{\cal
F} T'_B : B [G^o/X] \f B [G/X]$, donc $\a : O, y : B [G^o/X] \v_{\cal
F} (T'_B)y : B [G/X]$. D'o\`u $\a : O, y : B [G^o/X], x : B [G/X] \f C
[G/X] \v_{\cal F} (x) (T'_B)y : C [G/X]$. Or par hypoth\`ese
d'induction, $\a : O \v_{\cal F} T_C : C [G/X] \f C [G^o/X]$, donc $\a
: O, y : B [G^o/X], x : B [G/X] \f C [G/X] \v_{\cal F} (T_C)(x)
(T'_B)y : C [G^o/X]$. Par cons\'equent $\a : O \v_{\cal F} T_A = \l
x\l y (T_C) (x) (T'_B)y : A [G/X] \f A [G^o/X]$. La m\^eme
d\'emonstration se fait pour $\a : O \v_{\cal F} T'_A = \l x\l y
(T'_C) (x) (T_B)y : A [G^o/X] \f A [G/X]$.

-- Si $A = \q YB$, il faut d\'emontrer que $\a : O \v_{\cal F} \l x
(T_B)x : A [G/X] \f A [G^o/X]$ et $\a : O \v_{\cal F} \l x (T'_B)x : A
[G^o/X] \f A [G/X]$.  On a $x : \q YB [G/X] \v_{\cal F} x : B
[G/X]$. Comme par hypoth\`ese d'induction, $\a : O \v_{\cal F} T_B : B
[G/X] \f B [G^o/X]$, alors $\a : O, x : \q YB [G/X] \v_{\cal F} (T_B)x
: B [G^o/X]$. D'o\`u $\a : O, x : \q YB [G/X] \v_{\cal F} (T_B)x : \q
YB [G^o/X]$, et donc $\a : O \v_{\cal F} T_A = \l x (T_B)x : A [G/X]
\f A [G^o/X]$.  La m\^eme d\'emonstration se fait pour $\a : O
\v_{\cal F} T'_A = \l x (T'_B)x : A [G^o/X] \f A [G/X]$.  \cqfd \\
 
Notons que dans ce paragraphe, tous les $\l$-termes consid\'er\'es
sont typables et donc fortement normalisables.\\
   
{\bf D\'efinition} : On dit qu'un $\l$-terme est {\bf simple} s'il est
de la forme $(x)u_1...u_n$, avec $x$ une variable et $u_i$ ($1 \leq i
\leq n$) un $\l$-terme normal.\\
      
Dans le typage de $T$, on avait $\G \v_{\cal F} \d : A[G/X]$. On a le
lemme suivant :\\
 
{\bf Lemme 3.2} {\it On peut traiter seulement le cas o\`u $\d$ est un
$\l$-terme simple.}\\
  
{\bf Preuve :} Il faut \'etudier trois cas.

-- L'\'etape suivante dans le typage de $T$ est l'application de $\d$
\`a un terme $u$. Comme le terme $T$ est normal, $\d$ est un
$\l$-terme simple.

-- L'\'etape suivante dans le typage de $T$ est l'application de la
r\`egle $(\f_i)$, donc $\G- \{y : B\} \v_{\cal F} \l y \d : B \f
A[G/X]$, avec $y : B \in \G$.

Parall\`element on a, $\G- \{y : B\} \v_{\cal F} \l y (T'_A)\d : B \f
  A[G/X]$. Alors si $\a$ appartient \`a la forme normale de
  $(T'_A)\d$, il appartient \`a la forme normale de $\l y
  (T'_A)\d$. Comme le terme $T$ est normal, ce cas ne pose pas de
  probl\`eme.

 -- L'\'etape suivante dans le typage de $T$ est l'application d'un
  terme $u$ \`a $\d$, donc $\G' \v_{\cal F} u : A[G/X] \f B$ et $\G,
  \G' \v_{\cal F} (u)\d : B$.\\ Parall\`element on a, $\G, \G'
  \v_{\cal F} (u)(T'_A)\d : B$. $u$ est un $\l$-terme simple, car le
  terme $T$ est normal. Si $\a$ appartient \`a la forme normale de
  $(T'_A)\d$, alors il appartient \`a la forme normale de $(u)
  (T'_A)\d$. Par cons\'equent ce cas ne pose pas de probl\`eme. \cqfd
  \\

Dans la suite on suppose donc que $\d$ est un $\l$-terme simple.
 Montrons que $\a$ appartient \`a la forme normale de $(T'_A)\d$.\\
    
{\bf Lemme 3.3} {\it (i) Pour tout $\l$-terme simple $\D$, $\a$ est
libre dans la forme normale de $(T_A)\D$.

 (ii) Pour tout $\l$-terme simple $\D$, $\a$ est libre dans la forme
normale de $(T'_A)\D$.}\\
   
{\bf Preuve :} Par induction simultan\'ee sur $A$.

{\bf Preuve de (i) :} 

 -- Si $A = X$, alors $(T_X)\D \f_{\b} \l {\b} \l g (g)\D{\a}$, et
   donc c'est bon.

 -- Si $A = B \f C$, alors $(T_A)\D \f_{\b} \l y (T_C) (\D)
   (T'_B)y$. Si $X$ est libre dans $C$, alors, d'apr\`es l'hypoth\`ese
   d'induction, $\a$ est libre dans la forme normale de $(T_C) (\D)
   (T'_B)y$, donc $\a$ est libre dans la forme normale de
   $(T_A)\D$. Si $X$ n'est pas libre dans $C$, alors $(T_A)\D \f_{\b}
   \l y (\D) (T'_B)y$, et $X$ est libre dans $B$. Donc d'apr\`es (ii),
   $\a$ est libre dans la forme normale de $(T'_B)y$, d'o\`u $\a$ est
   libre dans la forme normale de $(T_A)\D$.

 -- Si $A = \q YB$, alors $(T_A)\D \f_{\b} (T_B)\D$, et donc par
   hypoth\`ese d'induction, on a le r\'esultat.

{\bf Preuve de (ii) :}

 -- Si $A = X$, alors $(T'_X)\D \f_{\b} (\D){\a}$ {\bf 1}.

 -- Pour les autres cas, on reprend la m\^eme preuve que (i). \cqfd \\
 
Dans la suite on va donner des cas particuliers qui prouvent que $E$
n'est pas un type sortie. \\
 
 {\bf Th\'eor\`eme 3.4} {\it Soient $t_1, ..., t_r$ des
$\l$-termes. Si $A$ se termine par $X$, alors la forme normale de
$(T'_A) \d t_1...t_r$ contient $\a$.}\\

  {\bf Preuve :} On a $A = \q \mbox{\boldmath$X_0$} (A_1 \f \q
\mbox{\boldmath$X_1$} (A_2\f...\f \q \mbox{\boldmath$X_{n-1}$} (A_n \f
\q \mbox{\boldmath$X_n$}X)...))$. \\ Donc $(T'_A)\d \f_{\b}
(T^{'}_{A_1 \f \q {\bf X_1} (A_2\f...\f \q {\bf X_{n-1}} (A_n \f \q
{\bf X_n}X)...)})\d \f_{\b}$\\ $\l y_1(T'_{\q {\bf X_1} (A_2\f...\f \q
{\bf X_{n-1}} (A_n \f \q {\bf X_n}X)...)}) (\d)(T_{A_1})y_1 \f_{\b}$
\\ $\l y_1(T'_{A_2\f...\f \q {\bf X_{n-1}} (A_n \f \q {\bf
X_n}X)...)}) (\d)(T_{A_1})y_1 \f_{\b}$\\ $\l y_1 \l y_2 (T'_{\q {\bf
X_2}(A_3\f...\f \q {\bf X_{n-1}} (A_n \f \q {\bf X_n}X)...)})
(\d)(T_{A_1})y_1 (T_{A_2})y_2 \f_{\b}$\\ $\l y_1...\l y_n
(T'_X)(\d)(T_{A_1})y_1...(T_{A_n})y_n$ {\bf 1} $\f_{\b}$ $\l y_1...\l
y_n (\d)(T_{A_1})y_1...(T_{A_n})y_n \a${\bf 1}.

 Dans la r\'eduction de $(T'_A) \d t_1...t_r$, trois cas peuvent se
produire :

 -- $r = n$, donc $(T'_A) \d t_1...t_r \f_{\b}
(\d)(T_{A_1})t_1...(T_{A_n})t_n\a$ {\bf 1}.

 -- $r < n$, donc $(T'_A) \d t_1...t_r \f_{\b} \l y_{r+1}...\l y_n
(\d)(T_{A_1})t_1...(T_{A_n})t_n\a$ {\bf 1}.

 -- $r > n$, donc $(T'_A) \d t_1...t_r \f_{\b}
(\d)(T_{A_1})t_1...(T_{A_n})t_n\a$ {\bf 1} $t_{n+1}...t_r$.
    
    On remarque que dans les trois cas, la forme normale
  de $(T'_A) \d t_1...t_r$ contient $\a$. \cqfd \\
  
Supposons donc que $A = \q \mbox{\boldmath$X_0$} (A_1 \f \q
\mbox{\boldmath$X_1$} (A_2\f...\f \q \mbox{\boldmath$X_{n-1}$} (A_n \f
\q \mbox{\boldmath$X_n$}Y)...))$, o\`u $Y$ est une variable
diff\'erente de $X$. On a le lemme suivant :\\
  
{\bf Lemme 3.5} {\it Soient $t_1, ..., t_r$ des $\l$-termes. Si l'un
des $A_i$ est \'egal \`a $X$, alors la forme normale de $(T'_A) \d
t_1...t_r$ contient $\a$.}\\
  
{\bf Preuve :} On a $(T'_A)\d \f_{\b}$

$ \l y_1...\l y_{n-1}(T'_{A_n \f \q {\bf X_n}Y})
(\d)(T_{A_1})y_1... (T_{A_i})y_i...(T_{A_{n-1}})y_{n-1} \f_{\b}$

 $\l y_1...\l y_{n-1}\l y_n (T'_{\q {\bf
X_n}Y})(\d)(T_{A_1})y_1...(T_{A_{n-1}})y_{n-1} (T_{A_n})y_n \f_{\b}$

$\l y_1...\l y_{n-1}\l y_n
(T'_Y)(\d)(T_{A_1})y_1...(T_{A_{n-1}})y_{n-1} (T_{A_n})y_n \f_{\b}$

$\l y_1...\l y_n (\d)(T_{A_1})y_1...(T_{A_n})y_n$. \\ Si $A_i = X$,
alors $(T_X) y_i \f_{\b} \l {\b} \l g (g)y_i{\a}$.\\ Donc $(T'_A)\d
\f_{\b} \l y_1...\l y_n (\d)(T_{A_1})y_1...\l {\b} \l g
(g)y_i{\a}...(T_{A_n})y_n$.\\ Comme $y_i$ reste en position
d'argument, alors par le m\^eme raisonnement du th\'eor\`eme 3.4, on
peut voir que la forme normale de $(T'_A) \d t_1...t_r$ contient
$\a$. \cqfd\\
  
Le lemme suivant est un raffinement du lemme 3.3.\\
   
{\bf Lemme 3.6} {\it Si $A$ se termine par $X$, alors la forme normale
de $(T_A)(x)u_1...u_n$ contient $\a$ qui n'est pas un argument de
$x$.}\\
  
{\bf Preuve :} Par induction sur $A$.

-- Si $A = X$, alors $(T_X)(x)u_1...u_n \f_{\b} \l {\b} \l g
((g)(x)u_1...u_n)\a$. Donc c'est bon.

-- Si $A = B \f C$, alors $C$ se termine par $X$, et
    $(T_A)(x)u_1...u_n$ $\f_{\b}$ 

$\l y (T_C)(x)u_1...u_n(T'_B)y$. Par hypoth\`ese d'induction sur $C$,
    on a le r\'esultat.

-- Si $A = \q YB$, alors $B$ se termine par $X$, et $(T_A)(x)u_1...u_n
\f_{\b} (T_B)(x)u_1...u_n$.  Donc, par hypoth\`ese d'induction sur
$B$, on a le r\'esultat.\cqfd\\
    
 {\bf Th\'eor\`eme 3.7} {\it Soient $t_1, ..., t_r$ des
$\l$-termes. Si l'un des $A_i$ ($1\leq i \leq n$) se termine par $X$,
alors la forme normale de $(T'_A) \d t_1...t_r$ contient $\a$.}\\
  
    {\bf Preuve :} On a $A = \q \mbox{\boldmath$X_0$} (A_1 \f \q
\mbox{\boldmath$X_1$}
   (A_2 \f...\f \q \mbox{\boldmath$X_{n-1}$}
    (A_n \f \q \mbox{\boldmath$X_n$}Y)...))$, avec $Y$ une variable.
Trois cas \`a examiner :

-- $A_i = X$, et on a donc le r\'esultat, d'apr\`es le lemme 3.5.

-- $A_i = C_i \f D_i$, et donc $D_i$ se termine par $X$. Par
cons\'equent :

  $(T'_A)\d \f_{\b} \l y_1...\l y_n (\d)(T_{A_1})y_1...
(T_{A_{i-1}})y_{i-1} (T_{C_i \f D_i})y_i...
  (T_{A_n})y_n \f_{\b}$\\ $\l y_1...\l y_n (\d)(T_{A_1})y_1...
(T_{A_{i-1}})y_{i-1} 
  \l z (T_{D_i})(y_i)(T'_{C_i})z...(T_{A_n})y_n$.

 Or d'apr\`es le lemme 3.6, la forme normale de
$(T_{D_i})(y_i)(T'_{C_i})z$ contient $\a$ qui n'est pas un argument de
$y_i$, donc la forme normale de $(T'_A) \d t_1...t_r$ contient $\a$.

-- $A_i = \q ZB_i$, et donc $B_i$ se termine par $X$.  On a $(T'_A)\d
\f_{\b} \\ \l y_1...\l y_n
(\d)(T_{A_1})y_1... (T_{B_i})y_i...(T_{A_n})y_n$.  D'o\`u, d'apr\`es
l'hypoth\`ese d'induction, la forme normale de $(T'_A) \d t_1...t_r$
contient $\a$. \cqfd \\

 {\bf D\'efinition :} On d\'efinit le syst\`eme ${\cal F}_F$ comme
\'etant le syst\`eme ${\cal F}$ o\`u on remplace la r\`egle $(\q_e)$
par la r\`egle :

\begin{center}
$(\q_e)_F \quad  \displaystyle\frac{ \G\v_{{\cal F}_F} t : \q XA}
{\G\v_{{\cal F}_F} t : A[G/X] }$
\end{center}

o\`u $A = \q \mbox{\boldmath$X_0$} (A_1 \f \q \mbox{\boldmath$X_1$}
   (A_2 \f...\f \q \mbox{\boldmath$X_{n-1}$} (A_n \f \q
   \mbox{\boldmath$X_n$}Y)...))$, avec $Y = X$ ou l'un des $A_i$ se
   termine par $X$.\\
    
    Alors on a le r\'esultat suivant :\\
    
  {\bf Th\'eor\`eme 3.8} {\it $E$ est un type sortie dans le syst\`eme
${\cal F}_F$ ssi $E$ est un type entr\'ee.}
  
 \section{Op\'erateurs de mise en m\'emoire}

{\bf D\'efinition :} Soient $D$ un type, et $T$ un $\l$-terme clos. On
dit que $T$ est un {\bf op\'erateur de mise en m\'emoire} (en
abr\'eg\'e {\bf o.m.m.}) {\bf pour $D$} ssi pour tout $\l$-terme $t$
tel que $\v_{\cal F} t : D$, il existe deux $\l$-termes $\t$, $\t'$,
avec $\t \simeq\sb{\b} \t'$ et $\v_{\cal F} \t' : D$, tel que pour
tout $\l$-terme $\th_t \simeq\sb{\b} t$, il existe une substitution
$\s$, telle que $(T)\th_t f \succ_f (f) \s (\t)$ o\`u $f$ est une
nouvelle variable.\\
  
{\bf D\'efinition :} Soit $\perp$ une constante de type
particuli\`ere. Pour toute formule $A$ de ${\cal F}$, on note par
$\neg A$ la formule $A \f \perp$, et par $A^*$ la formule obtenue en
rempla\c cant chaque formule atomique $R$ de $A$ par $\neg R$ ($A^*$
est dite la {\bf traduction de G\"odel de $A$}).\\
   
   On a le r\'esultat suivant (voir [8]) :\\
 
 {\bf Th\'eor\`eme 4.1} {\it Soit $D$ un type $\q^+$ tel que $\perp$
ne figure pas dans $D$. Si $\v_{\cal F} T : D^*\f ~\neg \neg D$, alors
$T$ est un o.m.m. pour $D$}.\\
    
Dans la preuve de ce th\'eor\`eme, on utilise deux propri\'et\'es
essentielles qui sont valables pour les types de donn\'ees
syntaxiques: l'une vient du fait que $D$ est un type entr\'ee et
l'autre est le th\'eor\`eme $2.1.1$ pour les types sorties.  D'o\`u le
r\'esultat suivant :\\
    
{\bf Th\'eor\`eme 4.2} {\it Soit $D$ un type de donn\'ees syntaxique
tel que $\perp$ ne figure pas dans $D$. Si $\v_{\cal F} T : D^* \f
\neg \neg D$, alors $T$ est un o.m.m. pour $D$.}\\

D'autre part on remarque que dans la preuve du th\'eor\`eme 4.2, on
utilise le type $D$ qui est \`a gauche de l'implication (dans
l'\'enonc\'e) comme type entr\'ee et le type \`a droite comme type
sortie. D'o\`u la d\'efinition et le r\'esultat suivants :\\
  
 {\bf D\'efinition :} On dit qu'un $\l$-terme clos $T$ est un {\bf
o.m.m. pour le couple de types $(E, S)$} ssi pour tout $\l$-terme $t$
v\'erifiant $\v_{\cal F} t : E$, il existe deux $\l$-termes $\t$,
$\t'$, avec $\t \simeq\sb{\b} \t'$ et $\v_{\cal F} \t' : S$, tel que
pour tout $\l$-terme $\th_t \simeq\sb{\b} t$, il existe une
substitution $\s$, telle que $(T)\th_t f \succ_f (f) \s (\t)$ o\`u $f$
est une nouvelle variable.\\
  
{\bf Th\'eor\`eme 4.3} {\it Soit $E$ un type entr\'ee et $S$ un type
sortie, tel que $\perp$ ne figure pas dans $S$. Si $\v_{\cal F} T :
E^* \f \neg \neg S$, alors $T$ est un o.m.m. pour $(E, S)$.}\\

{\bf Remarque :} La condition $S$ sortie est n\'ecessaire pour avoir
le th\'eor\`eme 4.3. En effet, si $S$ n'est pas sortie, alors il
existe un $\l$-terme normal $t$ contenant $x$ tel que $x : O \v_{\cal
F} t : S$. D'o\`u $x : E^* \v_{\cal F} t : S$, car $O \not \in Fv
(S)$, et donc $x : E^* \v_{\cal F} \l y (y)t : \neg \neg S$. Ce qui
fait que $\v_{\cal F} T = \l x \l y (y)t : E^* \f \neg \neg
S$. D'autre part $T$ n'est pas un o.m.m. pour $(E, S)$, car $(T)\t f
\p_f (f)t [\t/x]$ pour tout $\t \simeq\sb{\b} t$ et $t [\t/x]$
contient la variable $x$.\\

\begin{center}
REMERCIEMENTS
\end{center}
Nous remercions R. David pour ses conseils et ses remarques.


\begin{thebibliography}{99}

\bibitem{K} S. Farkh,
{\it Types de donn\'ees en logique du second ordre}.  
Th\`ese de doctorat, Universit\'e de Savoie, France (1998).
 
\bibitem{K} S. Farkh et K. Nour,
R\'esultats de compl\'etudes pour les types $\q^+$ du syst\`eme ${\cal F}$.
{\it C. R. Acad. Sci. Paris} {\bf 326 - I} (1998) 275-279.
 
\bibitem{K} J.-Y. Girard, Y. Lafont and  P. Taylor, {\em Proofs and Types}.
Cambridge University Press (1986).
 
\bibitem{K} J.-L. Krivine.
{\it Lambda-calcul, types et mod\`eles}. Masson, Paris (1990).

\bibitem{K} J.-L. Krivine,
Classical Logic, Storage Operators and Second Order
Lambda-Calculs.
{\it Annals of Pure and Applied Logic} {\bf 68} (1994) 53-78.

\bibitem{K} J.-L. Krivine, 
Op\'erateurs de mise en m\'emoire et traduction de G\"odel.  
{\it Archive for Mathematical Logic} {\bf 30} (1990) 241-267.

\bibitem{K} K. Nour,
{\it Op\'erateurs de mise en m\'emoire en lambda-calcul pur et typ\'e}.
Th\`ese de doctorat, Universit\'e de Savoie, France (1993).

\bibitem{K} K. Nour,
Op\'erateurs de mise en m\'emoire et types $\q$-positifs.
{\it RAIRO - Inform. Th\'eor. Appl.} {\bf 30} (1996) 261-293.

\bibitem{K} K. Nour, 
 Les $I$-types du syst\`eme ${\cal F}$.
{\it RAIRO - Inform. Th\'eor. Appl.} (to appear).
 
\end{thebibliography}
\end{document}